\documentclass[a4paper,10pt]{article}
\usepackage[utf8]{inputenc}
\usepackage{graphicx,amsmath,bbm,mathrsfs,amssymb,psfrag}

\begin{document}
\author{Laure Gouba\footnote{Laure Gouba is currently a visiting scientist at the Abdus Salam International 
Centre for Theoretical Physics (ICTP). A biographical sketch of the author can be 
read at the last page of the paper. Emails: lgouba@ictp.it; laure.gouba@gmail.com; lrgouba@yahoo.fr}}

\title{African Doctorates in Mathematics from Burkina-Faso}

\maketitle

\begin{abstract}
A catalogue of African Doctorates in Mathematics has been compiled and published 
in 2007 by the late Professor Paulus Gerdes \cite{paulus}. In this paper, we revise 
and update  the list of mathematicians from Burkina-Faso  compiled in \cite{paulus}.
Starting from a short description of Burkina-Faso, we brieffly 
mention the education system and the mathematics programme in Burkina-Faso
from the pre-school level to university level.
A list of mathematicians native of Burkina Faso is given. 
\end{abstract}

\section{Introduction }
 
Burkina Faso is located in the heart of West Africa in the sahel region.
It is bordered to the East by Niger, to the North and West by Mali, and
to the South by the Ivory Coast, Ghana, Togo and Benin. Its area is 274,000 $\textrm{km}^2$.
Its climate is dry and divided into two seasons: the dry season (from mid-October to mid-May) and 
the raining season (from mid-May to mid-October). Its average temperature is between 30 and 35 degrees.
The official language is French; there are also about sixty national languages, 
the three main being the Moor\'e, Fulfuld\'e (or Fulani) and Dioula. 
The main foreign languages spoken in Burkina Faso are English, German and Arabic.
Its capital is Ouagadougou (about 1, 6 million inhabitants), the second largest city 
Bobo-Dioulasso (600,000 inhabitants) is the economic capital. 
The population of Burkina-Faso is estimated at about 17 million in 2014.

The country was colonized by France in the late 19th century. 
The Colony of Upper Volta was created in 1919. 
In 1932, the Upper Volta was suppressed and divided between the Ivory Coast, Mali and Niger. 
It was reestablished as a territorial entity in 1947 and became independent in 1960. 
The Upper Volta has been renamed  Burkina Faso the 4th of August 1984.

Burkina Faso is a country of agriculture and  stock raising.
The very large majority (80\%) lives off agriculture and raising livestock. The rest 
of the population is involved in various commercial and craft activities, private business, 
industry, or works in the public service for the State.

The modern economy of Burkina Faso depends upon mines, 
processing industries of raw material such as fruits, vegetables, cereals, meat 
and textile industries based on cotton. Burkina Faso is considered as the biggest 
producer of cotton in Africa.

The education system in Burkina Faso has three sectors:
\begin{itemize}
 \item the formal sector, which includes pre-school, primary education, secondary education, higher 
 education and professional training organized in the context of schools;
 \item the non-formal sector which includes rural education and adult literacy. 
 This form of education is organized out of school context.
 \item Finally, the informal sector that takes into account the education received in the family circle or in a group.
\end{itemize}
The enrollment rate in 2014 was $72\%$ in primary school, $22\%$ in secondary school, and $6\%$ in higher education. 
These $6\%$, corresponding to a rate of about $0.3\%$ of the population, which is well below the standard value of
Unesco that is $2\%$.

We are interested to the formal sector of education. The preschool is optional and concern the range 
from 3 years old to 6 years old. The primary education starts from 6-7 years old to 12-13 years olds  
and lasts for six years. At the end of six years, 
the elementary student is required to pass the exam of the Certificate of Primary and Elementary 
Education (CEPE) which is required to access the secondary school.
The secondary education is divided into two steps. The first step of the formal education lasts four years.
At the fourth year the students pass the exam the first step studies (BEPC). The second step of the 
secondary school lasts 3 years and is divided in three mains options that are the general option, the technical 
option and the professional option. The last year of the second step of secondary school, the students pass 
an exam called BAC considered as the first degree of University which is compulsary for the registration at university.

Higher Education in Burkina-Faso began its structure in 1960s, 
with the country's accession to independence by establishing institutes or higher education schools.
In Burkina-Faso there are 4 national universities, 8 private universities and about 
65 Private higher education institutions.
The oldest university is the University of Ouagadougou established in 1974 and 
renamed in December 2015 as Universit\'e Ouaga I Pr Joseph Ki-Zerbo.

The paper is organized as follows. In section \ref{sec1}, we brieffly mention the mathematics programmes 
in Burkina starting with the general education system in Burkina-Faso. Then in section \ref{sec2}, we list the 
African mathematicians native of Burkina-Faso. Most of them are faculty members in Burkina-Faso and the rest 
are faculty members in universities abroad. 

\section{Mathematics in Burkina-Faso}\label{sec1}

 Long before the colonial penetration, mathematics existed in Burkina Faso, although they were not in a formal way. 
They were applied in everyday life. for example the counting system in dialects, the construction of traditional houses as round huts, 
decorations of traditional houses in geometric figures, making masks, etc ... Then there has been some evolution 
with the introduction of the western-style education system. Since then the mathematics enjoy a certain prestige in Burkina-Faso. 
The teaching style  and the programmes were quite similar to the ones in France until 1980s.
Since then, different approaches, initiated in Africa and supported by France have proposed 
slightly different programmes.

Let's brieffly mention mathematics programmes in the formal education system, the details being given 
in \cite{laure}.
\begin{itemize}
\item 
Mathematics begin at the preschool/kindergarden where the pupil learns to count, 
do activities and learn skills that will be extended in future mathematical 
learning like in primary school  with the introduction of  calculations,
subtraction, division. In Burkina Faso, the teaching of geometry is effective earlier at level of kindergarten. 
At this stage, the child must be able to recognize and classify certain objects by referring to their shapes.
\item The arithmetic, geometry and the metric system are the sub-disciplines of mathematics taught in primary school.
\item At the first step of secondary school that last 4 years, 
the programmes focus on the basics of pure geometry, arithmetic, vector calculus and algebra.
\item At the second step of secondary school that last 3 years the programmes of mathematics depend 
on the options that varies from scientific studies, literature studies or technical studies. For the scientific option, the programmes
of the first year of this level are based on numerical function 
of a real variable, equations and inequations in $\mathbb{R}$, vectors of the plan, geometry in the plane, 
transformation in plane, geometry in space and statistics. In the second year for the scientific option, the 
programs are based on algebraic and numerical problems, numerical sequences, numerical functions, angles and 
trigonometry, transformations in plane, geometry in space, statistics, enumerations. In the last year the programmes 
of mathematics are about arithmetic, probabilities, complex numbers, numerical sequences, numerical functions, integral 
calculation, planar curves, vector calculation and configuration (plane and space), transformations and configurations.
\item
Before 1974, year of the creation of the University of Ouagadougou,
and later the creation of the Institute of Mathematics and Physics (IMP), the first students in mathematics were trained abroad, 
mainly in France. Although some had begun their studies of mathematics in Africa, particularly in Dakar, 
Abidjan where there were the biggest african universities in 1960s -1970s, the first mathematicians completed their thesis in France.
In Burkina Faso, only the university of Ouagadougou has a complete training in mathematics. 
The current programme is Licence-Master-Doctorate programme labelled (LMD) started in 2009. 
More details about the LMD programme  in mathematics are given in \cite{laure}.
\end{itemize}
There are 4 research laboratories in mathematics:
\begin{itemize}
 \item Le laboratoire d'analyse num\'erique, d'informatique et de biomath\'ematiques (LANIBIO)\\
 Director: Professor Blaise Som\'e
 \item Le laboratoire de math\'ematiques et informatique (LAMI)\\
 Director: Professor Hamidou Tour\'e
 \item Le laboratoire de th\'eorie des nombres, alg\`ebre, g\'eom\'etrie alg\'ebrique, topologie alg\'ebrique (TN-AGATA)\\
 Director: Professor Moussa Ouattara
 \item Le Laboratoire de math\'ematiques, Universit\'e Polytechnique de Bobo-Dioulasso\\
 Director: Professor Marie Yves Th\'eodore Tabsoba
\end{itemize}

\section {African mathematicians native of Burkina-Faso}\label{sec2}

\begin{enumerate}
\item
Albert OUEDRAOGO  (Male)\\
1969: Doctorat de 3\`eme cycle\\
Title: Probl\`eme inverse de la diffusion et g\'en\'eralisation de l'\'equation de Marchenko 
( Inverse problem of diffusion and generalisation of Marchenko's equation).\\
Universit\'e Pierre et Marie Curie (Paris 6), France.\\
Directeur: J. L. Destouches\\
1981: Doctorat d'Etat\\
Title: Contr\^oles ponctuels de syst\`emes elliptiques et paraboliques d'ordre 2m:
application \`a un syst\`eme parabolique avec masse de Dirac (Punctual controls of elliptic 
and parabolic systems of order 2m).\\
Universit\'e Pierre et Marie Curie (Paris 6), France.\\
Directeur: Jacques Louis Lions
 \item
Akry KOULIBALY (Male) [Deceased]\\
1976 : Doctorat de 3\`eme cycle\\
Title: Alg\`ebres de Malcev de basses dimensions ( Malcev agebras of low dimensions).\\
Universit\'e de Montpellier 2, France.\\
Directeur: Artibano Micali\\
1984: Doctorat d'Etat\\
Title: Contributions \`a la th\'eorie de Malcev (Contributions to the theory of Malcev algebra).\\
Universit\'e Montpellier 2, France.\\
Directeur: Artibano Micali
\item 
Ousseynou NAKOULIMA (Male)\\
1977 : Doctorat de 3\`eme cycle\\
Title: Etude d'une in\'equation variationnelle bilat\'erale et d'un syst\`eme d'in\'equations  
quasi-variationnelles unilat\'erales associ\'ees.\\
Universit\'e de Bordeaux 1, France.\\
1981: Doctorat Sc. Math.\\
Title: In\'equations variationnelles et in\'equations quasivariationnelles bilat\'erales associ\'ees 
\`a des probl\`emes de jeux stochastiques \`a somme nulle ou non nulle, Universit\'e de Bordeaux 1, France.
\item
G\'erard KIENTEGA (Male)\\
1980: Doctorat de 3\`eme cycle\\                                                                        
Title: Sur les corps alg\'ebriques du quatri\`eme d\'egr\'e.\\
Universit\'e Pierre et Marie Curie (Paris  6), France. \\                                                                                                  
Directeur: Pierre Barrucand\\
1992:  PhD\\
Title: M\'etriques g\'en\'eralis\'ees et alg\`ebres  affinement compl\`etes.\\
Universit\'e de Montr\'eal, Canada.\\
Directeur: Ivo Rosenberg
\item 
Aboubakary Seynou (Male)\\
1981: Doctorat de 3\`eme cycle\\
Title: Compatibilit\'e temporelle,  Universit\'e Louis Pasteur de Strasbourg, France.\\
Directeur : Claude Dellacherie
\item 
Alfred TOURE (Male)\\
1981: Doctorat de 3\`eme cycle\\
Title: Divers aspects des connexions conformes.\\
Universit\'e Pierre et Marie Curie (Paris 6), France.\\
Directeur: Jacqueline Lelong Ferrand\\
1993: Doctorat Unique\\
Title: Geom\'etrie diff\'erentielle de certains fibr\'es unitaires.\\
Universit\'e de Montpellier 2, France.\\
Directeur: Jacques Lafontaine
\item
Hamidou TOURE (Male)\\
1982 : Doctorat de 3\`eme cycle\\
Title: Sur l'\'equation g\'en\'erale par la th\'eorie des semi-groupes 
non lin\'eaires dans L1 
(Non linear semi-group theory in L1 for a general equation).\\
Universit\'e de Franche Comt\'e, Besan\c{c}on, France.\\
Directeur: Philippe Benilan\\
1994 : Doctorat Unique\\
Title: Etude de probl\`emes fortement d\'eg\'ener\'es en une dimension d'espace 
(Study of strong degenerated parabolic problems  in one space dimension).\\
Universit\'e de Franche Comt\'e, Besan\c{c}on, France.\\
Directeur: Philippe Benilan\\
1995: Doctorat d'Etat\\
Title: Etude de probl\`emes paraboliques hyperboliques non lin\'eaires ( 
On nonlinear hyperbolic parabolic problems).\\
Universit\'e de Ouagadougou, Burkina-Faso.\\
Directeurs: Philippe Benilan, Albert Ouedraogo
\item
Dembo GADIAGA (Male)\\
1982: Doctorat de 3\`eme cycle\\ 
Title: Sur une classe de tests qui contient le test $\chi^2$:  le cas 
d'un processus stationnaire 
(On a class of tests which contain the $\chi^2$-test for a stationary process).\\
Universit\'e de Lille 1, France.\\
Directeur: Denis Bosq\\
2003: Doctorat d'Etat\\
Title: Test fonctionnel d'ajustement et de non influence pour des 
variables al\'eatoires d\'ependantes
(Functional tests  and no effects hypothesis for dependent  random variables).\\
Universit\'e de Ouagadougou, Burkina-Faso.\\
Directeurs: Denis Bosq, Albert Ouedraogo
\item 
Sabeko Marcel BONKIAN (Male)\\
1983: Doctorat de 3\`eme cycle\\
Title: Contribution \`a l'\'etude des mesures al\'eatoires du second ordre 
(Contribution to the study of random  measures of the second order).\\
Universit\'e des Sciences et Techniques de Lille 1, France.\\
 Directeur: Pierre Jacob
 \item 
Ba Amidou Boubacar YOBI (Male)\\
1983: Doctorat de 3\`eme cycle\\
Title: Contribution \`a l'\'etude des diagrammes de  De Finetti 
( Contribution to the study of De Finetti Diagramms ). \\
Universit\'e des Sciences et Techniques du Languedoc,
Montpellier, France.\\
Directeur: Artibano Micali
\item 
Blaise SOME (Male)\\
1984: Doctorat de 3\`eme cycle\\
Title: Identification, contr\^ole optimal et optimisation dans les syst\`emes 
diff\'erentiels compartimentaux 
(Identification, optimal control and optimisation in compartimental 
differential equation system).\\
Universit\'e Pierre et Marie Curie (Paris 6), France.\\
Directeur: Yves Cherruault\\
1994 : Doctorat d'Etat\\
Title: Algorithmiques num\'eriques et r\'esolution de probl\`emes de contr\^ole 
optimal et d'\'equations int\'egrales 
(Numerical algorithm and resolution of optimal control and integral equation problems).\\
Universit\'e de Ouagadougou, Burkina-Faso.\\
Directeur: Yves Cherruault
\item 
Longin SOME (Male)\\
1984: Doctorat de 3\`eme cycle\\
Title: Mise en oeuvre informatique de quelques m\'ethodes multigrilles dans le cadre de la m\'ethode 
des \'el\'ements finis (Computation of some multigrid methods under finite elements method).\\
Universit\'e Pierre et Marie Curie (Paris 6), France.\\
Directeur: Pierre Arnaud Raviard\\
2007: Doctorat Unique\\
Title: M\'ethode de grille mobile sous la m\'ethode des lignes pour la r\'esolution num\'erique d'\'equations 
aux d\'eriv\'ees partielles mod\'elisant des ph\'enom\`enes \'evolutifs.\\
Universit\'e de Ouagadougou, 
Burkina-Faso et  Facult\'e Polytechnique de Mons, Belgique.\\
Directeurs: Albert OUEDRAOGO, Philippe SAUCEZ  
\item
Marie Yves Th\'eodore TABSOBA (Male)\\
1987: Doctorat de 3\`eme cycle\\
Title: Complexit\'e de suites automatiques, Universit\'e Aix Marseille 2, France.\\
 1999:  Doctorat d'Etat\\
Title: Contribution \`a l'\'etude des suites automatiques, Universit\'e de Ouagadougou, Burkina-Faso.\\
Directeurs : G\'erard Rauzy et Akry Koulibaly
\item
Moussa OUATTARA (Male)\\
1988:  Doctorat Unique\\
Title: Alg\`ebre de Jordan et alg\`ebres g\'en\'etiques (Jordan algebras and genetic algebras).\\
Universit\'e de Montpellier 2, France. \\                         
Directeur: Artibano Micali\\
1991: Doctorat d'Etat\\
Title: Alg\`ebres de la g\'en\'etique des populations (Algebras of the genetics of population).\\
Universit\'e de Ouagadougou, Burkina-Faso.\\
 Directeurs: Artibano Micali et Akry Koulibaly
\item 
Kaka Bernard BONZI (Male)\\                                                        
1990: Doctorat Unique\\
Title: Etude des \'equilibres thermiques d'un supraconducteur, existence et stabilit\'e. 
Universit\'e de Nancy , France.\\
Directeur: Lanchon-Ducauquis H\'el\`ene
\item
Kalifa TRAORE (Male)\\
1990: Doctorat 3\`eme cycle\\
Title: Cohomologie des alg\`ebres de Malcev (Cohomology of Malcev algebra).\\
Universit\'e de Ouagadougou, Burkina-Faso.\\
Directeur: Akry Koulibaly\\
2006: PhD\\
Title: Etudes des pratiques math\'ematiques d\'evelopp\'ees  en contexte par 
les Siamous du Burkina Faso.\\
Universit\'e du Quebec \`a Montr\'eal, Canada.\\
Directeurs: Nadine Bednarz,  Philippe Jonnaert.
\item
Bourama TONI (Male)\\
1994 : PhD\\
Title: Bifurcations de p\'eriodes critiques locales.\\
Universit\'e de Montr\'eal, Canada.\\
Directeur: Christiane Rousseau.
\item 
Sado TRAORE (Male)\\
1994: Doctorat Unique\\
Title: Approche variationnelle de la dualit\'e quasi convexe (Variational 
approach to quasi convex duality).\\
Universit\'e d'Avignon et des pays du  Vaucluse, Avignon,  France.\\
Directeur: Michelle Volle.
\item 
Nakelgbamba Boukary PILABRE (Male)\\
1995: Doctorat 3\`eme cycle\\
Title: Sur la Lie admissibilit\'e de la dupliqu\'ee non commutative  d'une alg\`ebre 
(On the Lie admissibility of the noncommutative  duplication of an algebra).\\
Universit\'e de Ouagadougou, Burkina Faso.\\
Directeur:  Akry Koulibaly\\
2011 : Doctorat de l'Universit\'e de Ouagadougou.\\
Title: Dupliqu\'ee et quelques structures alg\'ebriques. \\
Universit\'e de Ouagadougou, Burkina-Faso.\\
Directeur: Moussa Ouattara
\item
C\^ome Jean Antoine BERE (Male)\\
1997: Doctorat de 3\`eme cycle\\
Title: Superalg\`ebres de Malcev.\\
Universit\'e de Ouagadougou, Burkina Faso.\\
Directeur: Akry Koulibaly
\item
Pierre Clovis NITIEMA (Male)\\
1998: PhD\\
Title: Approximation des puissances sectionn\'ees et des classes de fonctions \`a l'aide 
de polyn\^omes alg\'ebriques.\\
Universit\'e d'Etat de Dniepropetrovsk, Facult\'e de M\'ecanique et de Math\'ematiques, Ukraine.\\
Directeur: Motornyi Vitaly Pavlovitch.
\item
Mamadou SANGO (Male)\\
1998:  Ph D (Th\`ese Unique)\\
Title: Valeurs propres et vecteurs propres de probl\`emes elliptiques nonautoadjoints 
avec un poids ind\'efini.\\
Universit\'e de Valenciennes et du Hainaut-Cambr\'esis, France.\\
Directeur  : Serge Nicaise 
\item
Sibiri TRAORE (Male)\\
1998: PhD\\
Title: Conception d'une architecture neuronique de grande taille inspir\'ee de 
l'anatomie et de la physiologie 
du cerveau humain  pour la commande de robots.\\
Universit\'e Laval, Canada.\\
Directeur: Michel Guillot
\item
Joseph BAYARA (Male)\\
1999: Doctorat 3\`eme cycle\\ 
Title: Sur les alg\`ebres d'\'evolution (On evolution  algebras).\\
Universit\'e de Ouagadougou, Burkina-Faso.\\
Directeur : Moussa Ouattara\\
2013: Doctorat de l'Universit\'e de Ouagadougou\\
Title: Alg\`ebres train et d\'erivations dans les alg\`ebres associatives.\\
Universit\'e de Ouagadougou, Burkina-Faso.\\
Directeur: Moussa Ouattara
\item 
Marie Fran\c{c}oise OUEDRAOGO (Female)\\
1999: Doctorat 3\`eme cycle\\
Title: Sur les superalg\`ebres triples de Lie (On Lie Triple superalgebras).\\
Universit\'e de Ouagadougou, Burkina-Faso.\\
Directeur:  Akry Koulibaly\\
2009: Doctorat de l'Universit\'e Blaise Pascal de Clermont-Ferrand\\
Title: Extension of the canonical trace and associated determinants.\\
Universit\'e Blaise Pascal de Clermont-Ferrand, France.\\
Directeurs: Sylvie Paycha et Akry Coulibaly
\item 
Andr\'e  CONSEIBO (Male)\\
2001: Doctorat 3\`eme cycle \\
Title: Alg\`ebre de Berstein  monog\`ene et \'equation  diff\'erentielle (Monogenic Berstein 
algebra and differential equation).\\
Universit\'e de Ouagadougou, Burkina-Faso.\\
Directeur: Moussa Ouattara\\
2013: Doctorat Unique\\
Title: Alg\`ebres train de degr\'e 2 et d'exposant 3 (Train algebras of degree 2 and exponent 3).\\
Universit\'e de Ouagadougou, Burkina-Faso.\\
Directeur: Moussa Ouattara
\item 
Stanislas OUARO (Male)\\
2001: Doctorat Unique \\
Title:  Etude de probl\`emes  elliptiques-paraboliques non lin\'eaires en une dimension d'espace 
 (On non-linear  elliptic parabolic problems in one space dimension).\\
 Universit\'e de Ouagadougou, Burkina-Faso.\\
Directeur:  Hamidou Tour\'e
\item 
Lucie Patricia ZOUNGRANA (Female)\\
2001: Doctorat Unique\\
Title: Sur les sous boucles de Cartan d'une boucle homog\`ene et les sous alg\`ebres 
de Cartan d'une alg\`ebre triple de Lie [On Cartan  subloops of homogenous loop 
and Cartan subalgebra of a Lie triple algebra].\\
Universit\'e de Ouagadougou, Burkina-Faso.\\
Directeur: Akry Koulibaly
\item 
Mahamadi Jacob WARMA (Male)\\
2002: PhD\\
Title: The Laplacian with Robin boundary conditions on arbitrary domains.\\
University of Ulm, Germany.\\
Directeur : Wolfgang Arendt.
\item 
Oumar TRAORE (Male)\\
2002: Doctorat Unique\\
Title: Contr\^ole de probl\`emes dynamiques de la population (Control in population dynamics problems).\\
Universit\'e de Ouagadougou, Burkina-Faso.\\
Directeur: Albert Ouedraogo.
\item 
Idrissa KABORE (Male)\\
2004: Doctorat Unique\\
Title: Combinatoire des mots et caract\'erisation de certaines classes.\\
Universit\'e de Ouagadougou, Burkina-Faso.\\
Directeur: Marie Yves Th\'eodore Tapsoba.
\item 
Somdouda SAWADOGO (Male)\\
2005: Doctorat unique\\
Title: Contr\^olabilit\'e de syst\`emes dynamiques \`a deux temps. 
Application \`a la th\'eorie des sentinelles
( Controlability of dissipative systems. Application to the sentinel theory). \\
Universit\'e de Ouagadougou, Burkina-Faso.\\
Directeurs:  Ousseynou Nakoulima et Albert Ouedraogo
\item 
Ousseni SO (Male)\\
2005: Doctorat Unique\\
Title: Mod\'elisation Math\'ematique et simulation num\'erique de la physiologie et du contr\^ole optimal 
des \'echanges gazeux respiratoires chez l'homme  
(Mathematics modeling and numerical simulation of physiology and optimal control of 
the human respiratory system).\\
Universit\'e de Ouagadougou, Burkina-Faso.\\
Directeur: Blaise Som\'e
\item 
Laure GOUBA (Female)\\
2005: PhD \\
Title: Th\'eories de Jauge Ab\'eliennes Scalaires et Spinorielles \`a 1+1 Dimensions: 
une Etude non Perturbative.\\
Institut de Math\'ematiques et de Sciences Physiques (IMSP), Porto-Novo, 
Universit\'e d'Abomey Calavi (UAC), R\'epublique du B\'enin.\\
Directeurs: Jan Govaerts,  Norbert Mahouton Hounkonnou
\item 
Balira Ousmane  KONFE (Male)\\
2005: Th\`ese Unique\\
Title: Optimisation globale: M\'ethode Alienor.\\
Universit\'e de Ouagadougou, Burkina-Faso.\\
Directeur: Blaise Som\'e
\item 
Genevi\`eve BARRO (Female)\\
2005: Th\`ese Unique\\
Title: Contribution \`a la r\'esolution num\'erique de quelques probl\`emes de r\'eactions 
diffusion non lin\'eaires.\\
Universit\'e de Ouagadougou, Burkina-Faso.\\
Directeurs:  Benjamin Mampassi, Blaise Som\'e.
\item 
Mikailou COMPAORE (Male)\\
2007: Doctorat Unique \\
Title: Sur les d\'eformations infinit\'esimales des flots et la g\'eom\'etrie diff\'erentielle des fibres unitaires de 
certains espaces sym\'etriques de rang 1.\\
Universit\'e de Ouagadougou, Burkina-Faso.\\
Directeur: Edmond Fedida
\item 
 W. Jacob YOUGBARE (Male)\\
2007: Th\`ese Unique \\
Title: Data envelopment analysis/ M\'ethodologie, th\'eorie et relation avec l'optimisation multicrit\`eres: 
application aux syst\`emes de l'enseignement de base, de la sant\'e et de quelques entreprises du Burkina-Faso.\\
Universit\'e de Ouagadougou, Burkina-Faso.\\
Directeur: Blaise Som\'e
\item 
Jean de Dieu ZABSONRE (Male)\\
2008: Doctorat Unique\\
Title: Mod\`eles visqueux en s\'edimentation et stratification ( obtention formelle, stabilit\'e th\'eorique 
et sch\'emas volumes finis bien \'equilibr\'es).\\
Universit\'e de Savoie, France.\\
Directeurs: Hamidou Tour\'e,  Didier Bresch,  E. Fernandez-Nieto.
\item 
Aboudramane GUIRO (Male)\\
2009: Doctorat Unique\\
Title: Sur quelques probl\`emes d'observateurs: Applications \`a certains mod\`eles d' \'ecosyst\`eme aquatique.\\
Universit\'e de Ouagadougou, Burkina-Faso.\\
Directeurs: Hamidou Tour\'e et  Abderrahman IGGIDR.
\item 
Adama OUEDRAOGO (Male)\\
2009: Doctorat Unique\\
Title: Solutions renormalis\'ees pour des probl\`emes paraboliques fortement d\'eg\'en\'er\'es: cas isotrope et non isotrope.\\
Universit\'e de Ouagadougou, Burkina-Faso.\\
Directeurs: Hamidou Tour\'e et Mohamed Maliki
\item 
Issa ZABSONRE (Male)\\
2009: Doctorat Unique\\
Title: Contributions \`a l'\'etude d'une classe d'\'equations int\'egro-diff\'erentielles \`a retard en dimension infinie. \\
Universit\'e de Ouagadougou, Burkina-Faso.\\
Directeurs: Hamidou Tour\'e et  Khalil Ezzinbi.
\item 
Gilbert BAYILI (Male)\\
2009: Doctorat Unique\\
Title: Contr\^ole  des coefficients de singularit\'es et contr\^olabilit\'e exacte dans un domaine polygonal avec fissures.\\
Universit\'e de Ouagadougou, Burkina-Faso.\\
Directeurs: Hamidou Tour\'e et  Mary Teuw Niane.
\item 
Pascal ZONGO (Male)\\
2009: Th\`ese Unique \\
Title: Mod\'elisation math\'ematique de la dynamique de transmission du paludisme.\\
Universit\'e de Ouagadougou, Burkina-Faso.\\
Directeur:  Blaise Som\'e
\item 
Seydou Eric TRAORE (Male)\\
2010: Doctorat Unique\\
Title: Ensemble et Syst\`emes Flous et Applications \`a Quelques Mod\`eles de D\'ecisions en Sciences Environnementales.\\
Universit\'e de Ouagadougou, Burkina-Faso.\\
Directeurs: Hamidou Tour\'e et  Akry Koulibaly
\item 
Hermann Wendpayande Baudoin SORE (Male)\\ 
2010: Doctorat de l'Universit\'e Hambourg\\
Title: The Dold-Kan Correspondance and coalgebras structures.\\ 
Universitat Hamburg, Allemagne.\\
Directeur: Birgit Richter.
\item 
 Youssouf PARE (Male)\\
2010: Th\`ese Unique\\
Title: R\'esolution de quelques \'equations fonctionnelles par la m\'ethode num\'erique SBA.\\
Universit\'e de Ouagadougou, Burkina-Faso.\\
Directeur: Blaise Som\'e
\item 
 Elis\'ee GOUBA (Male)\\
2010: Th\`ese Unique\\
Title: Identification de param\`etres dans les syst\`emes distribu\'es \`a donn\'ees manquantes.
Mod\`eles math\'ematiques de la performance en sport.\\
Universit\'e de Ouagadougou, Burkina-Faso et Universit\'e des Antilles et de la Guyane, France.\\
Directeurs: Olivier HUE, Ousseynou Nakoulima et Blaise Som\'e
\item 
Diakarya  BARRO (Male)\\
2010: Th\`ese Unique \\
Title: Contributions \`a la mod\'elisation statistique des valeurs extr\^emes multivari\'ees.\\
Laboratoire LANIBIO, UFR-SEA, Universit\'e de Ouagadougou, Burkina-Faso.\\
Directeurs:  Dossou-Gb\'et\'e Simplice et Blaise Som\'e
\item 
Safimba SOMA (Male)\\
2011: Doctorat Unique \\
Title: Etude de probl\`emes elliptiques non lin\'eaires avec des donn\'ees mesures. \\
Universit\'e de Ouagadougou, Burkina-Faso.\\
Directeurs:  Stanislas Ouaro et Nouredine Igbida
\item 
Blaise KONE (Male)\\
2011: Doctorat Unique\\
Title: Etudes de probl\`emes anisotropiques non lin\'eaires.\\
Universit\'e de Ouagadougou, Burkina-Faso.\\
Directeur : Stanislas Ouaro
\item 
Isma\"el NYANKINI (Male)\\
2012: Doctorat Unique\\
Title: Etude de probl\`emes elliptiques non lin\'eaires sous des conditions assez g\'en\'erales sur les donn\'ees.\\ 
Universit\'e de Ouagadougou, Burkina-Faso.\\
Directeur:  Stanislas Ouaro
\item 
Olivier  SAWADOGO (Male)\\
2012: Th\`ese Unique\\
Title: Mod\'elisation hydrog\'eologique : Ecoulement en milieux poreux, fractur\'es, 
Probl\`eme inverse et transport de polluants.\\
Universit\'e de Ouagadougou, Burkina-Faso.\\
Directeur:  Blaise Som\'e
\item 
 Malicki ZOROM (Male)\\
2012: Th\`ese Unique\\
Title: Mod\'elisation compartimentale: dynamique de la vuln\'erabilit\'e socio-\'economique des ruraux sah\'eliens  \`a la 
variabilit\'e climatique et contr\^ole optimal d'un mod\`ele de type metapopulation du paludisme.\\
Universit\'e de Ouagadougou, Burkina-Faso.\\
Directeur: Blaise Som\'e
\item 
Boureima SANGARE (Male)\\
2012: Th\'ese Unique\\
Title: Adaptation dynamique de maillage pour la r\'esolution des \'equations aux d\'eriv\'ees partielles \'evolutives.\\ 
Universit\'e de Ouagadougou, Burkina-Faso et Universit\'e des Sciences, des Techniques et des Technologies de Bamako, Mali.\\
Directeurs : Ouateni Diallo et Longin Som\'e
\item 
Moumini KERE (Male) \\
2012: Th\`ese Unique \\
Title: Contr\^olabilit\'e simultan\'ee et application \`a un probl\`eme d'identification simultan\'ee de param\`etre.\\
Universit\'e de Ouagadougou, Burkina-Faso.\\
Directeurs: Ousseynou Nakoulima et Blaise Som\'e
\item 
Sadou TAO (Male)\\
2012: Th\`ese Unique\\
Title: Contr\^olabilit\'e de syst\`emes paraboliques et application aux sentinelles.\\
Universit\'e de Ouagadougou, Burkina-Faso.\\
Directeurs: Ousseynou Nakoulima et Blaise Som\'e
\item 
Victorien Fourtoua KONANE (Male)\\
2013: Doctorat Unique\\
Title: Etude de Syst\`emes Dynamiques (Mod\'elisation de batteries non rechargeables).\\
Universit\'e de Ouagadougou,  Burkina-Faso.\\
Directeurs: Dembo Gadiaga et Ingemar Kaj
\item 
Bila Adolphe KYELEM (Male)\\
2013: Doctorat Unique\\
Title: Contribution \`a l'\'etude d'existence de solutions p\'eriodiques pour une classe de probl\`emes 
d'\'evolution \`a retard et applications.\\
Universit\'e de Ouagadougou, Burkina-Faso.\\
Directeurs : Stanislas Ouaro et Khalil Ezzinbi
\item 
Ibrahim NONKANE (Male)\\
2013: Doctorat de l'Universit\'e de Ouagadougou\\
Title: Geom\'etrie des modules et des modules diff\'erentiels associ\'es aux repr\'esentations du groupe sym\'etrique.\\ 
Universit\'e de Ouagadougou, Burkina-Faso.\\
Directeurs:  Rikard Bogvad et G\'erard Kientega
\item 
Duni Yegbonoma Fr\'ed\'eric ZONGO (Male)\\
2013: Doctorat Unique\\
Title: Etude de probl\`emes anisotropiques nonlin\'eaires et  d'\'equations quasi relativistes de type Choquard. \\
Universit\'e de Ouagadougou, Burkina-Faso. \\
Directeurs: Stanislas Ouaro et  Michael Melgaard
\item 
Francis BASSONO (Male)\\
2013: Th\`ese Unique\\
Title: Etude de quelques \'equations fonctionnelles par les m\'ethodes SBA, d\'ecompositionnelle d'Adomian et des perturbations. \\  
Universit\'e de Ouagadougou, Burkina-Faso.\\
Directeurs:  Gabriel Bissanga et Blaise Som\'e
\item 
Dalomi BAHAN (Male)\\
2013: Th\`ese Unique\\
Title: M\'ethode Ali\'enor pour la r\'esolution des probl\`emes d'optimisation quadratiques en variables binaires: 
applications et complexit\'e. \\
Universit\'e de Ouagadougou, Burkina-Faso.\\
Directeur:  Blaise Som\'e
\item 
Kounhinir SOME (Male)\\
2013: Th\`ese Unique\\
Title: Nouvelle m\'eta heuristique bas\'ee sur la m\'ethode Ali\'enor pour la r\'esolution des probl\`emes 
d'optimisation multi objectif: th\'eorie et applications.\\
Universit\'e de Ouagadougou, Burkina-Faso.\\
Directeurs: Blaise Som\'e et Berthold Ulungu
\item 
Arouna OUEDRAOGO (Male)\\
2014: Doctorat Unique\\
Title: Etude de probl\`emes elliptiques et paraboliques non lin\'eaires.\\
Universit\'e de Ouagadougou, Burkina-Faso.\\
Directeur: Stanislas Ouaro
\item 
Salifou NIKIEMA (Male)\\
2015: Doctorat de l'Universit\'e de Ouagadougou\\
Title: Polyn\^omes \`a coefficients entiers : hauteurs et bornes pour les facteurs.\\
Universit\'e de Ouagadougou, Burkina-Faso. \\
Directeur: G\'erard Kientega
\end{enumerate}

{\bf Biographical Sketch}

Laure Gouba is native of Burkina-Faso where she was born and grew up. 
She is married to Riccardo Fantoni, an Italian physicist.

Laure studied mathematics at University of Ouagadougou and obtained  her 
master degree in applied Mathematics in 1999. She also taught Mathematics in 
High School as part time activity from October 1996 to May 2001.
She has been selected Full member of the Organization for Women in Science 
for the Developing World (OWSD) in 1999.
From 2001 till 2005 she was a PhD student at Institut de Math\'ematiques et de Sciences 
Physiques (IMSP)/Porto-Novo in Benin. In 2005 she visited UCL (Unversit\'e Catholique de Louvain
in Belgium) as an exchange student. She obtained her PhD in Theoretical/Mathematical Physics from the University 
of Abomey Calavi in November 2005. 

In 2006 she visited the African Institute for Mathematical Sciences (AIMS) in Muizenberg /South Africa
and the Abdus Salam International Centre for Theoretical Physics (ICTP) in Trieste/Italy.
From September 2006 till June 2008, she was a postdoctoral fellow and coordinator of teaching 
assistants at AIMS. She has been selected Associate member of ICTP for the period 2007 - 2012.
From July 2008  to September 2010, she was a postdoctoral fellow at the 
National Institute for Theoretical Physics  (NITheP) in  Stellenbosch/ South Africa.
Laure is member of the International Chair in Mathematical Physics and Applications (ICMPA), 
Cotonou /Benin.

Since October 2010 she is a visiting scientist at ICTP. 
Laure is pursuing her research in Theoretical/ Mathematical physics and 
she is also interested in science and education. She is author of about 20 
scientific papers and a referee for Journal of Mathematical Physics A and Journal of Physics A: Math-Theor.

\end{document}